\documentclass[11pt]{article}
\usepackage{amsfonts}
\usepackage{graphicx}
\usepackage{color}
\usepackage{amsmath,amssymb,comment}
\usepackage{verbatim}
\newtheorem{theo}{Theorem}[section]
\newtheorem{lem}[theo]{Lemma}

\newtheorem{prop}[theo]{Proposition}

\newcommand{\mysection}[1]{\section{#1} \setcounter{equation}{0}}
    \makeatletter
\def\@fnsymbol#1{\ensuremath{\ifcase#1\or *\or \ddagger\or
   \mathsection\or \mathparagraph\or \|\or **\or \dagger\dagger
   \or \ddagger\ddagger \else\@ctrerr\fi}}
    \makeatother
\newcommand{\proof}{{\sc Proof.} \quad}
\newcommand{\proofc}{{\sc Proof} \ }
\newcommand{\be}{\begin{equation} \label}
\newcommand{\ee}{\end{equation}}
\newcommand{\bea}{\begin{eqnarray}\label}
\newcommand{\eea}{\end{eqnarray}}
\newcommand{\bas}{\begin{eqnarray*}}
\newcommand{\eas}{\end{eqnarray*}}
\newcommand{\bit}{\begin{itemize}}
\newcommand{\eit}{\end{itemize}}
\newcommand{\qed}{\hfill$\Box$ \vskip.2cm}
\newcommand{\nn}{\nonumber}
\newcommand{\R}{\mathbb{R}}
\newcommand{\N}{\mathbb{N}}
\newcommand{\pO}{\partial\Omega}

\newcommand{\wto}{\rightharpoonup}
\newcommand{\wsto}{\stackrel{\star}{\rightharpoonup}}

\newcommand{\io}{\int_\Omega}
\newcommand{\na}{\nabla}
\newcommand{\Del}{\Delta}
\newcommand{\del}{\delta}
\newcommand{\al}{\alpha}

\newcommand{\sig}{\sigma}
\newcommand{\pa}{\partial}
\newcommand{\bom}{\overline{\Omega}}
\newcommand{\Om}{\Omega}

\newcommand{\hs}{\hspace*}

\newcommand{\vp}{\varphi}
\newcommand{\lbal}{\left\{ \begin{array}{l}}
\newcommand{\lball}{\left\{ \begin{array}{ll}}
\newcommand{\ear}{\end{array} \right.}

\newcommand{\abs}{\\[5pt]}

\newcommand{\tm}{T_{max}(\al,\bet)}
\renewcommand{\al}{\alpha}
\newcommand{\bet}{\beta}
\newcommand{\gam}{\gamma}
\newcommand{\uab}{u_{\al\bet}}
\newcommand{\vab}{v_{\al\bet}}
\newcommand{\wab}{w_{\al\bet}}
\newcommand{\zab}{z_{\al\bet}}
\newcommand{\ug}{u_\gam}
\newcommand{\vg}{v_\gam}
\newcommand{\wg}{w_\gam}
\newcommand{\zg}{z_\gam}
\newcommand{\sab}{S_{\al\bet}}
\newcommand{\tab}{T_{\al\bet}}
\newcommand{\ow}{\overline{w}}
\newcommand{\uinf}{u_\infty}
\newcommand{\vinf}{v_\infty}

\begin{document}
\enlargethispage{10mm}
\title{Phenotype switching in chemotaxis aggregation models controls the spontaneous emergence of large densities}
\author
{
Kevin J Painter\footnote{kevin.painter@polito.it}\\
{\small Politecnico di Torino, DIST}\\
{\small 10125 Torino, Italy}
\and
Michael Winkler\footnote{michael.winkler@math.uni-paderborn.de}\\
{\small Universit\"at Paderborn, Institut f\"ur Mathematik,}\\
{\small 33098 Paderborn, Germany} }
\date{}
\maketitle
\begin{abstract}
\noindent

We consider a phenotype-switching chemotaxis model for aggregation, in which a chemotactic population is capable of switching back and forth between a chemotaxing state (performing chemotactic movement) and a secreting state (producing the attractant). We show that the switching rate provides a powerful mechanism for controlling the densities of spontaneously emerging aggregates. Specifically, in two- and three-dimensional settings it is shown that when both switching rates coincide and are suitably large, then the densities of both the chemotaxing and the secreting population will exceed any prescribed level at some points in the considered domain. This is complemented by two results asserting the absence of such aggregation phenomena in corresponding scenarios in which one of the switching rates remains within some bounded interval.\abs
\end{abstract}
\newpage
\section{Introduction}\label{intro}
The Keller-Segel system for chemotactic movement \cite{KS70} constitutes a classical model of mathematical biology, noted for its capacity to exhibit self-organisation. When formulated with two variables -- a chemotactic population and its diffusive chemoattractant -- and making the assumption that the population secretes its own attractant, positive feedback is inserted into this system that allows clusters to form through self-reinforcement. Studies of this system form a significant literature, ranging from biological and sociological applications of the model \cite{painter} to studies into the intricate analytical properties \cite{hillen_painter2009,bellomo}. In terms of the latter, resolving the question of whether solutions to the system exist globally-in-time or blow up has been of particular interest.

This simple model implicitly assumes homogeneity across the population, i.e. that all members have equivalent capacity to perform chemotaxis and/or produce attractant. These assumptions of homogeneity, however, can lie at odds with the heterogeneity observed in natural systems, where a population will often span a spectrum of {\em phenotypes}, such as variable chemotactic sensitivity in the context of chemotactic behaviour \cite{salek}. The presence of phenotypic heterogeneity may confer numerous advantages, such as adaptability within fluctuating environments, but must also be considered within the wider context of energy balances \cite{keegstra}: performing one activity (e.g. movement) has an inevitable energy cost, meaning less resources may be available for other activities (e.g. protein synthesis and transport), i.e. there exists a {\em trade-off}. In light of this, recent models have been developed that extend the Keller and Segel framework to include phenotypic heterogeneity and trade-offs, either through a population that spans a continuous spectrum of phenotypes \cite{lorenzi} or a population subdivided into discrete phenotypic states \cite{macfarlane}. The latter considered a self-organisation type scenario, specifically focusing on a trade-off in which the population either performed chemotaxis ({\em chemotaxing state}) or produced the attractant ({\em secreting state}), but could not perform both simultaneously. Linear stability analysis of this model demonstrated that pattern formation is possible, but that the rate of switching between the two states becomes a crucial factor in the dynamics. Numerical analyses confirmed these predictions and, moreover, hinted that solutions exist globally-in-time, yet a formal analysis was not attempted.

Motivated by this, we will focus on a family of problems based on the system introduced in \cite{macfarlane}. In the following, $\uab$ represents the chemotaxing phenotype, $\wab$ represents the secreting phenotype and $\vab$ is the chemoattractant. Specifically, we shall consider the problems

\be{0ab}
	\left\{ \begin{array}{ll}	
	\pa_t \uab = \Del \uab - \na\cdot (\uab \na \vab) - \al \uab + \bet \wab,
	\qquad & x\in\Om, \ t>0, \\[1mm]
	\pa_t \vab = \Del \vab - \vab + \wab,
	\qquad & x\in\Om, \ t>0, \\[1mm]
	\pa_t \wab = \Del \wab - \bet \wab + \al \uab,
	\qquad & x\in\Om, \ t>0, \\[1mm]
	\frac{\pa\uab}{\pa\nu}=\frac{\pa\vab}{\pa\nu}=\frac{\pa\wab}{\pa\nu}=0,
	\qquad & x\in\pO, \ t>0, \\[1mm]
	\uab(x,0)=u_0(x), \quad \vab(x,0)=v_0(x), \quad \wab(x,0)=w_0(x),
	\qquad & x\in\Om,
	\end{array} \right.
\ee
with $\al>0$ and $\bet>0$ representing the rate of switching from chemotaxing to secreting and secreting to chemotaxing, respectively.\abs 

Viewed from a mathematical perspective, (\ref{0ab}) differs from the corresponding initial-boundary value problem for the classical Keller-Segel system
\be{KS}
	\lbal
	u_t = \Del u - \chi\na\cdot (u\na v), \\[1mm]
	v_t = \Del v - v + \mu u,
	\ear
\ee
by including a form of indirectness with respect to the signal production process.
Whereas the classical system is well-known for its ability to generate unbounded densities within the chemotaxing population
(\cite{herrero_velazquez}, \cite{win_JMPA}), the inclusion of indirectness within the signal dynamics has been found to significantly reduce a trend toward singular destabilisation. In particular, in \cite{fujie_senba} the presence of additional diffusion-induced dissipation of the form in (\ref{0ab}) 
was shown to suppress any blow-up phenomena in a closely related problem in three (and lower)-dimensional settings. Meanwhile, a substantial collection of additional studies has provided comparable findings for a number of nearby extensions
(\cite{dong_peng}, \cite{liu_wang}, \cite{ren_liu}, \cite{wu}, \cite{xing}, \cite{ye_wang}, \cite{zhang_niu_liu}). Lacking, as far as we are aware, is a rigorous understanding into the extent to which regularised systems of this nature retain core features of taxis-driven destabilisation, despite the absence of unboundedness phenomena.\abs
{\bf Main results.}
The purpose of the present manuscript is to confirm that the coefficients $\al$ and $\beta$ of the zero
order contributions to (\ref{0ab}) play a key role in this regard. In particular, simultaneously increasing these parameters can result in the spontaneous emergence of structures which, while not singular, do lead to population densities of arbitrary size.\abs
As a preliminary step, in Section \ref{sect2} we confirm that, within frameworks of physically relevant dimension, blow-up does not occur in any of the populations in (\ref{0ab}), in the sense that for any choice of the parameters therein, and for all reasonably regular nonnegative initial data,
a global bounded classical solution can always be found:
\begin{prop}\label{prop1}
  Let $n\le 3$ and $\Om\subset\R^n$ be a bounded domain with smooth boundary, let $\al>0$ and $\bet>0$, and suppose that
  \be{init}
	u_0\in W^{1,\infty}(\Om),
	v_0\in W^{1,\infty}(\Om)
	\mbox{ and } 
	w_0\in W^{1,\infty}(\Om)
	\quad \mbox{are nonnegative}.
  \ee
  Then there exist uniquely determined nonnegative functions
  \be{1.1}
	\lbal
	\uab\in C^0(\bom\times [0,\infty)) \cap C^{2,1}(\bom\times (0,\infty)), \\[1mm]
	\vab\in \bigcap_{q>n} C^0([0,\infty);W^{1,q}(\Om)) \cap C^{2,1}(\bom\times (0,\infty))
	\qquad \mbox{and} \\[1mm]
	\wab\in C^0(\bom\times [0,\infty)) \cap C^{2,1}(\bom\times (0,\infty))
	\ear
  \ee
  such that $(\uab,\vab,\wab)$ forms a classical solution of (\ref{0ab}).
  Moreover, this solution is bounded in the sense that for each $q>n$ there exists $C(q)>0$ such that
  \be{1.2}
	\|\uab(\cdot,t)\|_{L^\infty(\Om)}
	+ \|\vab(\cdot,t)\|_{W^{1,q}(\Om)}
	+ \|\wab(\cdot,t)\|_{L^\infty(\Om)}
	\le C(q)
	\qquad \mbox{for all } t>0.
  \ee
\end{prop}
In the light of a corresponding result for a relative of (\ref{0ab}), obtained in \cite{fujie_senba}, this
observation is not in itself particularly surprising. The approach taken to derive this, however, is designed in such a way that it can be artlessly further developed so as to reveal that the bounds in (\ref{1.2}) are actually uniform within large parts of the parameter region $(0,\infty)^2$.\abs
In fact, as a first statement in this direction, we
show that provided $\al$ is located in a fixed interval, then the inequality in (\ref{1.2}) remains stable for arbitrarily large $\bet$:
\begin{prop}\label{prop41}
  Let $n\le 3$, and assume (\ref{init}).
  Then for each $q>n$ and any $\al^\star>0$ and $\del>0$, there exists $C(q,\al^\star,\del)>0$ such that 
  \bea{41.1}
	& & \hs{-20mm}
	\|\uab(\cdot,t)\|_{L^\infty(\Om)}
	+ \|\vab(\cdot,t)\|_{W^{1,q}(\Om)}
	+ \|\wab(\cdot,t)\|_{L^\infty(\Om)}
	\le C(q,\al^\star,\del) \nn\\[1mm]
	& & \hs{30mm}
	\qquad \mbox{for all $t>0$ and any choice of $\al\in (0,\al^\star]$ and $\bet\ge \del$.}
  \eea
\end{prop}
Likewise, under the corresponding scenario that $\bet$ is located in a fixed interval, unboundedness phenomena can also be ruled out in the limit $\al\to\infty$:
\begin{prop}\label{prop42}
  Let $n\le 3$, and suppose that (\ref{init}) holds.
  Then whenever $q>n, \bet^\star>0$ and $\del\in (0,\bet^\star)$, one can find $C(q,\bet^\star,\del)>0$ with the property that 
  \bea{42.1}
	& & \hs{-20mm}
	\|\uab(\cdot,t)\|_{L^\infty(\Om)}
	+ \|\vab(\cdot,t)\|_{W^{1,q}(\Om)}
	+ \|\wab(\cdot,t)\|_{L^\infty(\Om)}
	\le C(q,\bet^\star,\del) \nn\\[1mm]
	& & \hs{30mm}
	\qquad \mbox{for all $t>0$, any $\al>0$ and each } \bet\in [\del,\bet^\star].
  \eea
\end{prop}
In stark contrast to the above results, however, our final result will indicate that when {\em both} parameters in (\ref{0ab})
simultaneously diverge, then arbitrarily large population densities may spontaneously emerge in two- and three-dimensional
balls. 
To address this in as simple as possible setting, we shall concentrate here on the prototypical case when the large numbers $\al$ and $\bet$ precisely coincide, and hence
consider the one-parameter sub-family of (\ref{0ab}) given by
\be{0}
	\left\{ \begin{array}{ll}	
	\pa_t \ug = \Del \ug - \na\cdot (\ug \na \vg) - \gam \ug + \gam \wg,
	\qquad & x\in\Om, \ t>0, \\[1mm]
	\pa_t \vg = \Del \vg - \vab + \wg,
	\qquad & x\in\Om, \ t>0, \\[1mm]
	\pa_t \wg = \Del \wg - \gam \wg + \gam \ug,
	\qquad & x\in\Om, \ t>0, \\[1mm]
	\frac{\pa\ug}{\pa\nu}=\frac{\pa\vg}{\pa\nu}=\frac{\pa\wg}{\pa\nu}=0,
	\qquad & x\in\pO, \ t>0, \\[1mm]
	\ug(x,0)=u_0(x), \quad \vg(x,0)=v_0(x), \quad \wg(x,0)=w_0(x),
	\qquad & x\in\Om,
	\end{array} \right.
\ee
with $\gam>0$.\abs
Within this context, for some fixed initial data the corresponding solutions can be seen to undergo a quite drastic unboundedness phenomenon, as to be confirmed in Section \ref{sect3}:
\begin{theo}\label{theo12}
  Let $n\in \{2,3\}, R>0$ and $\Om=B_R(0)\subset\R^n$.
  Then there exist $T>0$ and radially symmetric initial data which satisfy (\ref{init}), and which are such that whenever
  $(\gam_j)_{j\in\N} \subset (0,\infty)$ is unbounded, one can find a subsequence $(\gam_{j_k})_{k\in\N}$ with the property
  that for the corresponding solutions of (\ref{0}) from Proposition \ref{prop1} we have
  \be{12.1}
	\|u_{\gam_{j_k}}\|_{L^\infty(\Om\times (0,T))} \to \infty
	\qquad \mbox{as } k\to\infty
  \ee
  and
  \be{12.2}
	\|w_{\gam_{j_k}}\|_{L^\infty(\Om\times (0,T))} \to \infty
	\qquad \mbox{as } k\to\infty.
  \ee
\end{theo}
We remark here that the taxis-driven spontaneous emergence of large densities has been discovered within extreme parameter constellations for several Keller-Segel type systems. However, unlike the majority of preceding literature in this regard
(cf., e.g.~\cite{lankeit_exceed}, \cite{kang_stevens}, \cite{fuest_heihoff}, \cite{WWX_CVPDE}, \cite{win_JNLS} or \cite{win_JDE2022}),
the observation that leads to Theorem \ref{theo12} does not refer to parameter limits which are singular in the sense that either fast or slow diffusion limits are involved; indeed, the underlying core mechanism here appears to be more complex, in that the dominance of chemotaxis over diffusion is enforced by an enhanced zero order interaction, despite the marked non-degeneracy of all migration processes.
\mysection{Numerical motivation}\label{numerics}
To motivate Propositions \ref{prop1}-\ref{prop42} and Theorem \ref{theo12} we compare numerical solutions of the switching model (\ref{0}) with those of the formulation
\be{0inf}
	\left\{ \begin{array}{ll}	
	\pa_t \uinf = \Del \uinf - \frac{1}{2} \na\cdot (\uinf \na \vinf),
	\qquad & x\in\Om, \ t>0, \\[1mm]
	\pa_t \vinf = \Del \vinf - \vinf + \frac{1}{2} \uinf,
	\qquad & x\in\Om, \ t>0, \\[1mm]
	\frac{\pa\uinf}{\pa\nu}=\frac{\pa\vinf}{\pa\nu}=0,
	\qquad & x\in\pO, \ t>0, \\[1mm]
	\uinf(x,0)=u_0(x), \quad \vinf(x,0)=v_0(x),
	\qquad & x\in\Om,
	\end{array} \right.
\ee
The above can be regarded as the classical companion to (\ref{0}), featuring a homogeneous population,  $\uinf$, that secretes its own chemoattractant, $\vinf$. A standard linear stability analysis of the uniform steady state solutions for (\ref{0inf}) predicts that self-organisation can occur, given sufficient population mass. However, subsequent dynamics vary according to the dimension, $n$: when $n=1$ solutions exist globally-in-time, yet when $n\ge 2$ finite-time blow-up can occur. Intuitively, (\ref{0inf}) could conceivably arise in the limit $\gamma \rightarrow \infty$ of (\ref{0}): switching between chemotaxis-only and secreting-only states becomes instantaneous, effectively generating a single-state population that simultaneously produces the attractant and performs chemotaxis; the half-factors reflect a division of time between the two states.  

\medskip
To substantiate this intuition we perform a numerical exploration for $n=1,2,3$. Generally, numerical studies in $n=2,3$ are compromised by computational cost, particularly if high accuracy is required. To circumvent this such that simulations in $n=1,2,3$ can be performed to an equivalent degree of resolution, each of the cases $n=2,3$ are restricted to an assumed radial symmetry. Specifically, for $n=1$ we consider the interval $\Omega = [0,R]$, while in $n=2,3$ we consider $\Omega = {\cal {B}}_R^n$ (the $n$-dimensional ball of radius $R$) and radially symmetric initial data, allowing reduction to the radial line $[0,R]$. All simulations utilise the initial distributions $\ug(x,0) = \wg(x,0) = \rho/2$, $\vg(x,0) = \rho(1+0.001e^{-x^2})/2$ for (\ref{0}) or $\uinf(x,0) = \rho$, $\vinf(x,0) = \rho(1+0.001e^{-x^2})/2$ for (\ref{0inf}); here, $x$ is used to refer to either the position along the interval $[0,R]$ in $n=1$, or the radial position in $n=2,3$. The parameter $\rho$ measures the population mass and linear stability analysis (see \cite{macfarlane}) can be used to determine the {\em autoaggregation} parameter spaces. Specifically, these are regions in $(R, \rho,\gamma)$-space for (\ref{0}) and $(R,\rho)$-space for (\ref{0inf}) in which the uniform steady state is unstable to the inhomogeneous perturbation and aggregations are predicted to emerge. Regarding the numerical implementation itself, we utilise the \texttt{pdepe} solver in {\sc{Matlab}}, which discretises in space to yield a system of ODEs to be integrated in time (using \texttt{ode15s}). The above initial conditions bias high density aggregate emergence to the origin, which can then be exploited by performing a spatial discretisation on a non-uniform mesh. Specifically we discretise $[0,R]$ into a mesh of $N+1$ points $X_0, X_1, \hdots, X_N$, where $X_i = R i^2/N^2$. This concentrates grid points near the origin, where high resolution is desirable for the potentially steep gradients forming at the point of mass concentration. We set $N=10^4$ and note that simulations with twice or half the number yielded quasi-identical results. Regardless of these measures, computation of $z(x,t)$ is still expected to fail in certain instances, for example due to emergence of extremely large densities; this is particular to be expected in (\ref{0inf}) for $n=2,3$, where finite time blow up is possible. Any such scenarios of numerical failure are classified as ``numerical blow-up'' and we compute up to the critical time $t_c$ at which the numerical scheme fails.

\begin{figure}[t!]
\centering
\includegraphics[width=\textwidth]{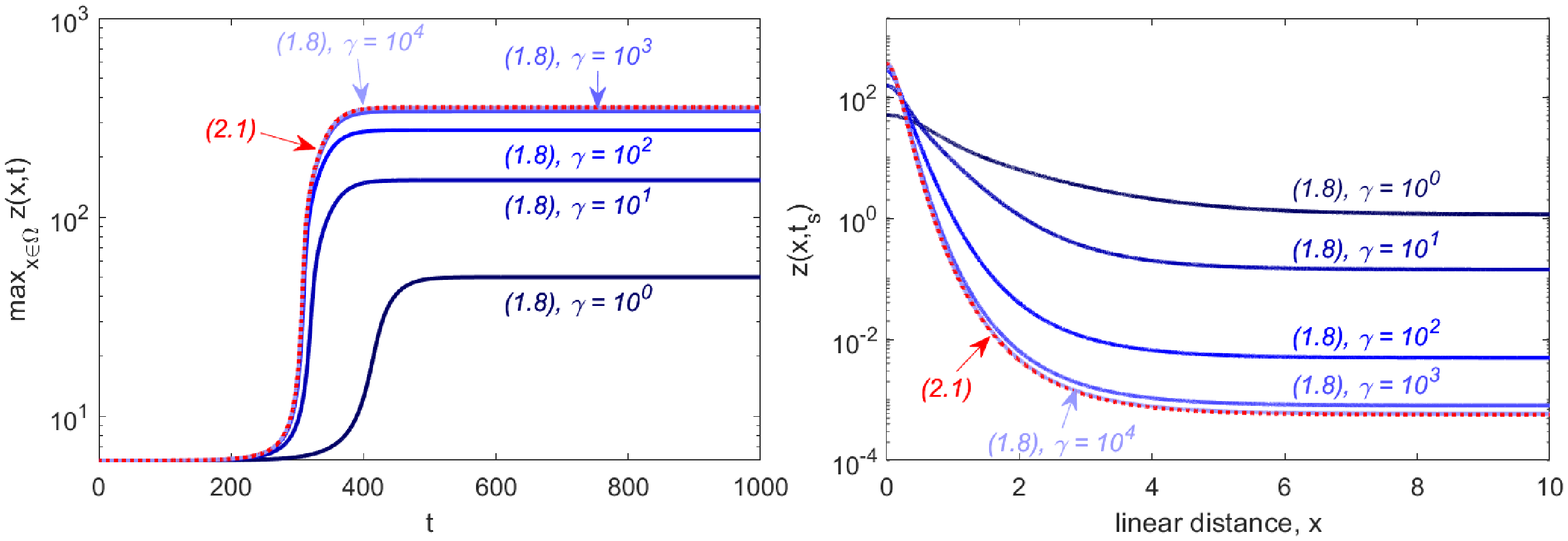}
\caption{{\bf Numerical solutions for (\ref{0}) and (\ref{0inf}) when $n=1$.} Left. Time-evolution of the maximum total population density, $\max_{x\in\Omega} z(x,t)$. Dotted red line shows the solution to (\ref{0inf}), while blue shaded solid lines indicate solutions to (\ref{0}) with different choices for $\gamma$. Right. Total population density distribution at $t = t_s$, i.e. $z(x,t_s)$, where $t_s$ is a point in time at which there is negligible change to the solution (numerical steady state). For these simulations we set $R=10$ and $\rho = 6$, with initial conditions as stated in the text.}\label{figure1d}
\end{figure}

\medskip
We first consider the case $n=1$, see Figure \ref{figure1d}. Note that instead of plotting individual variables, we consider the total population density $z(x,t)$, i.e. $z_{\gamma} = u_{\gamma}+w_{\gamma}$ for (\ref{0}) and $z_{\infty} = u_{\infty}$ for (\ref{0inf}). When $n=1$ we expect global existence for (\ref{0inf}). Selecting parameters from the autoaggregation parameter space for (\ref{0inf}), we correspondingly observe the growth of a cluster (Figure \ref{figure1d}, dotted-red line) at the  boundary $x=0$ and converging towards a high density aggregate; the dotted red-line in Figure \ref{figure1d} (right) shows the long-term solution to (\ref{0inf}), computed at a point when negligible solution change is observed. Similar dynamics are observed for the switching model (\ref{0}), with the various solid blue-shaded lines in Figure \ref{figure1d} corresponding to distinct choices for $\gamma$. For smaller $\gamma$ the aggregate grows over a slower timescale and settles towards a nonuniform steady state with lower peak density, this delay/reduction resulting from the cost of switching between chemotaxis-only and producing-only states. Note that when $\gamma$ is reduced below a critical threshold, autoaggregation is not possible and solutions rather evolve towards the uniform steady state (see \cite{macfarlane} for details). As $\gamma$ is increased, however, faster growth of the cluster occurs and we observe convergence in both space and time between the numerical solutions to (\ref{0}) and those of (\ref{0inf}). This provides substance to our intuition that (\ref{0}) converges to (\ref{0inf}) as $\gamma \rightarrow \infty$. 

\begin{figure}[t!]
\centering
\includegraphics[width=\textwidth]{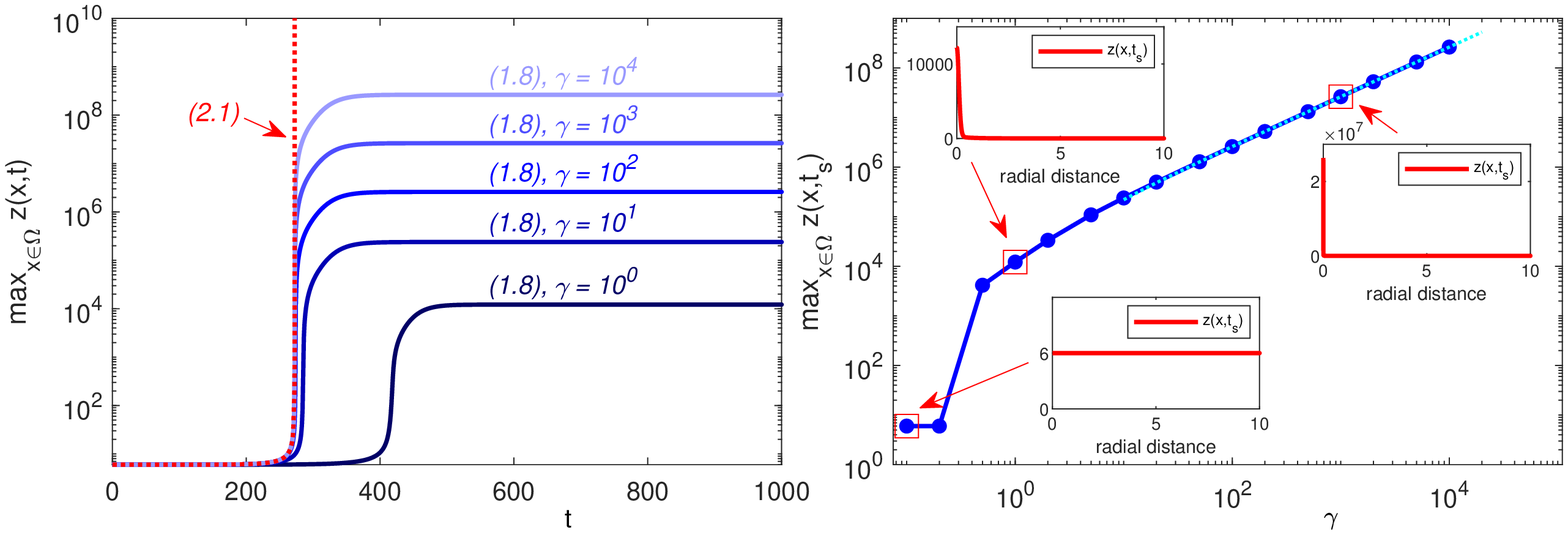}
\caption{{\bf Numerical solutions for (\ref{0}) and (\ref{0inf}) when $n=2$, assuming radial symmetry.} Left. Time-evolution of the maximum total population density, $\max_{x\in\Omega} z(x,t)$. Dotted red line shows the solution to (\ref{0inf}), while blue shaded solid lines indicate solutions to (\ref{0}) using different choices for $\gamma$. Note that the numerical solution of (\ref{0inf}) fails at $t_c \sim 272.269$ (numerical blow-up). Right. Maximum of the total population density distribution at $t = t_s$, i.e. $\max_{x\in\Omega} z(x,t_s)$, where $t_s$ is a point in time at which there is negligible change to the solution (numerical steady state). The insets show the total population density distribution at the same time point, $z(x,t_s)$, for selected $\gamma$. Note further that the straight line fit through the data for values $\gamma \ge 10^1$ is shown via the dotted cyan line. For these simulations we set $L=10$ and $\rho = 6$, with initial conditions as stated in the text.}\label{figure2d}
\end{figure}

\medskip
We next explore the $n=2$ case. Standard theory predicts that solutions to (\ref{0inf}) will blow-up in finite-time given sufficient initial mass. Correspondingly, upon simulation of (\ref{0inf}) we observe the initial growth of a cluster that rapidly accelerates into a highly concentrated peak: the growth of this peak is indicated by the red dotted line in Figure \ref{figure2d} (left). Numerical blow-up occurs at $t_c \sim 272.269$, beyond which no computation is possible within the limits of the numerical scheme invoked. In stark contrast, corresponding solutions to (\ref{0}) appear to exist globally in-time and no numerical blow-up was observed across the range of $\gamma$ used (five orders of magnitude). Tracking the maximum density of the cluster, we observe a similar initial trend in which the cluster first increases slowly before entering a phase of rapidly accelerating growth. Notably, increasing $\gamma$ leads to a convergence between solutions to (\ref{0}) and solutions to (\ref{0inf}) prior to $t=t_c$, with the time of most rapid acceleration in cluster growth for  (\ref{0}) coinciding with the numerical blow-up time observed for (\ref{0inf}). Beyond this point, however, the growth rate of the clusters that form in (\ref{0}) is curtailed, and solutions begin to converge to what appear to be bounded nonuniform steady state solutions. To quantify this, we plot the maximum density of the steady state cluster distribution that forms for (\ref{0}) as a function of $\gamma$, Figure \ref{figure2d} (right). For larger $\gamma$ we observe an almost perfectly linear relationship between the size of $\gamma$ and the maximum density, implying that there is no clear upper bound to the maximum density of solutions to (\ref{0}) if $\gamma$ is allowed to become arbitrarily high. 

\medskip
As a final control, we repeat the numerical simulations for the case $n=3$. The results largely mirror those observed for $n=2$, where we observe that while numerical blow-up occurs for the system (\ref{0inf}) (here, for $t_c \sim 276.756$), solutions to (\ref{0}) instead evolve to a tightly aggregated solution with a bounded upper density. As for $n=2$, increasing $\gamma$ leads to a convergence between solutions to the model 
(\ref{0}) and the solutions to (\ref{0inf}) prior to the time of numerical blow up. Note that the solution densities formed in $n=3$ are notably higher than those observed in $n=2$ for the corresponding value of $\gamma$.  Despite this, however, the numerical solver was still able to compute solutions for (\ref{0}) without numerical blow-up, with growth of the cluster slowing and the solution stabilising into a tightly aggregated steady state solution.

\begin{figure}[t!]
\centering
\includegraphics[width=\textwidth]{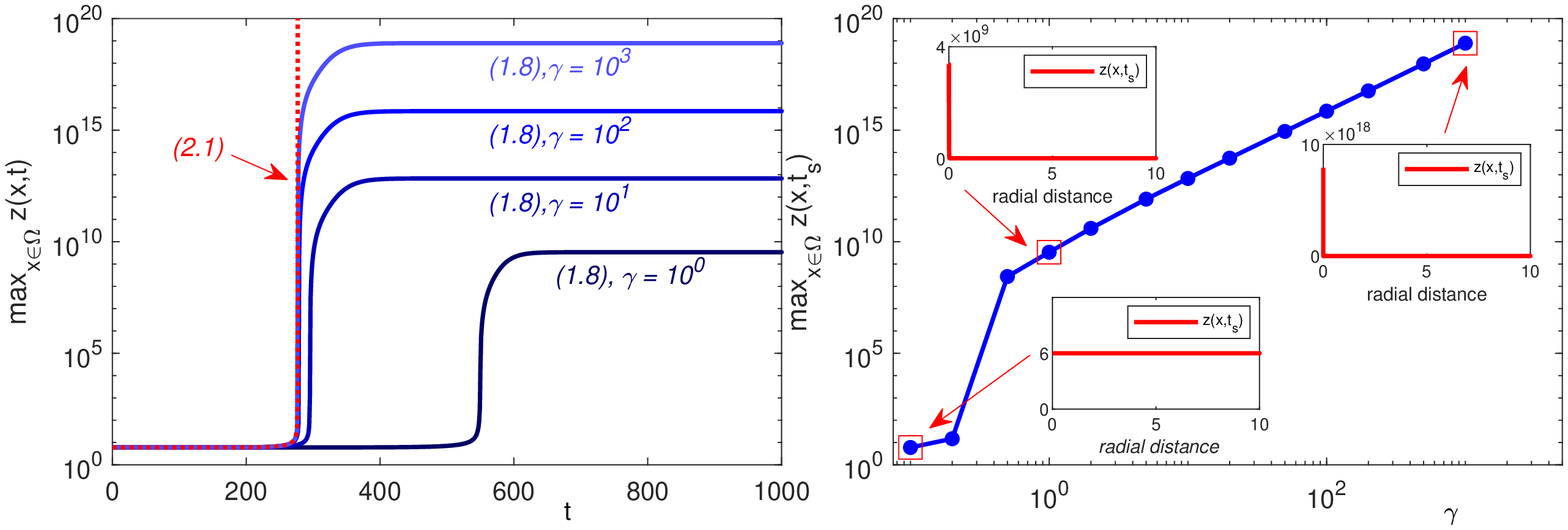}
\caption{{\bf Numerical solutions for (\ref{0}) and (\ref{0inf}) when $n=3$, assuming radial symmetry.} Left. Time-evolution of the maximum total population density, $\max_{x\in\Omega} z(x,t)$. Dotted red line shows the solution to (\ref{0inf}), while blue shaded solid lines indicate solutions to (\ref{0}) using different choices for $\gamma$. Note that the numerical solution of (\ref{0inf}) fails at $t_c \sim 276.757$ (numerical blow-up). Right. Maximum of the total population density distribution at $t = t_s$, i.e. $\max_{x\in\Omega} z(x,t_s)$, where $t_s$ is a point in time at which there is negligible change to the solution (numerical steady state). The insets show the total population density distribution at the same time point, $z(x,t_s)$, for selected $\gamma$. For these simulations we set $R=10$ and $\rho = 6$, with initial conditions as stated in the text.}\label{figure3d}
\end{figure}

\mysection{Boundedness results. Proofs of Propositions \ref{prop1}, \ref{prop41} and \ref{prop42}}\label{sect2}
\subsection{Local existence. Basic properties of $\uab+\wab$}
Let us first recall some standard theory on short-term solvability in triangular taxis-type parabolic systems (\cite{amann},
\cite{horstmann_win}) to state, without further comment, the following result on local existence and extensibility of solutions to the general problem (\ref{0ab}), along with some simple observations concerned with the sum of the corresponding two
population densities addressed therein.
\begin{lem}\label{lem_loc}
  Let $n\ge 1$, assume (\ref{init}), and let $\al>0$ and $\bet>0$.
  Then there exist $\tm\in (0,\infty]$ and a classical solution $(\uab,\vab,\wab)$ of (\ref{0ab}) in $\Om\times (0,\tm)$,
  uniquely determined by the inclusions
  \bas
	\lbal
	\uab \in C^0(\bom\times (0,\tm)) \cap C^{2,1}(\bom\times (0,\tm)), \\[1mm]
	\vab \in \bigcap_{q>n} C^0([0,\tm);W^{1,q}(\Om)) \cap C^{2,1}(\bom\times (0,\tm))
	\qquad \mbox{and} \\[1mm]
	\wab \in C^0(\bom\times (0,\tm)) \cap C^{2,1}(\bom\times (0,\tm)),
	\ear
  \eas
  such that $\uab\ge 0, \vab\ge 0$ and $\wab\ge 0$ in $\Om\times (0,\tm)$, and that
  \bea{ext}
	& & \hs{-10mm}
	\mbox{if $\tm<\infty$, \qquad then \qquad} \\[1mm]
	& & \hs{6mm}
	\limsup_{t\nearrow\tm} \Big\{
	\|u(\cdot,t)\|_{L^\infty(\Om)}
	+ \|v(\cdot,t)\|_{W^{1,q}(\Om)}
	+ \|w(\cdot,t)\|_{L^\infty(\Om)}
	\Big\}
	=\infty
	\quad \mbox{for all $q>n$.}
  \eea
  Moreover,
  \be{z}
	\zab(x,t):=\uab(x,t)+\wab(x,t),
	\qquad x\in\bom, \ t\in [0,\tm),
  \ee
  satisfies
  \be{0z}
	\left\{ \begin{array}{ll}	
	\pa_t \zab = \Del \zab - \na \cdot (\uab \na\vab),
	\qquad & x\in\Om, \ t\in (0,\tm), \\[1mm]
	\frac{\pa\zab}{\pa\nu}=0,
	\qquad & x\in\pO, \ t\in (0,\tm), \\[1mm]
	\zab(x,0)=u_0(x)+w_0(x),
	\qquad & x\in\Om,
	\end{array} \right.
  \ee
  and we have
  \be{mass}
	\io \zab(\cdot,t)=\io u_0 + \io w_0
	\qquad \mbox{for all } t\in (0,\tm).
  \ee
\end{lem}
Whenever $\al$ and $\bet$ as well as $(u_0,v_0,w_0)$ have been fixed, in what follows we shall (without explicitly stating) let $\tm$ as well as $\uab, \vab,\wab$ and $\zab$ be as obtained in Lemma \ref{lem_loc}.
For unambiguity in the notation, hereon we agree on choosing the norms in first order Sobolev spaces on $\Om$ according to the definition 
$\|\vp\|_{W^{1,q}(\Om)}:=\big\{ \io |\na \vp|^q + \io |\vp|^q\big\}^\frac{1}{q}$ for $\vp\in W^{1,q}(\Om)$ and $q\ge 1$.
\subsection{A general boundedness criterion. Estimates for individual solutions of (\ref{0ab})}
The core of the analysis in this section can now be found in the following result on boundedness in (\ref{0ab}),
conditional in the sense that, inter alia, an $L^1$ boundedness feature of the first solution components is presupposed.
Having this statement at hand will not only enable us to derive Proposition \ref{prop1} as a fairly direct consequence,
but also forms a key ingredient for the arguments on parameter-independent boundedness claimed in Propositions
\ref{prop41} and \ref{prop42}.
Our argument in this context is almost exclusively based on a suitably efficient exploitation of well-known
smoothing properties enjoyed by the Neumann heat semigroup on $\Om$; in fact, this will be seen to be possible on the mere basis
of the $L^1$ information in (\ref{31.2}), if the spatial dimension satisfies the $n\le 3$ assumption, explicitly relied on only in the following lemma.
\begin{lem}\label{lem31}
  Assume that $n\le 3$ and that (\ref{init}) holds. 
  Then for each $q>n$, $\del>0$ and $K>0$, one can find $C(q,K,\del)>0$ with the property that whenever $\al>0,\bet\ge\del$
  and $t_0\in [0,\tm)$ are such that
  \be{31.2}
	\al \|\uab(\cdot,t)\|_{L^1(\Om)} \le K
	\qquad \mbox{for all } t\in (t_0,\tm)
  \ee
  and
  \be{31.3}
	\|\uab(\cdot,t_0)\|_{L^\infty(\Om)}
	+ \|\vab(\cdot,t_0)\|_{W^{1,q}(\Om)}
	+ \|\wab(\cdot,t_0)\|_{L^\infty(\Om)}
	\le K,
  \ee
  then
  \be{31.4}
	\|\uab(\cdot,t)\|_{L^\infty(\Om)}
	+ \|\vab(\cdot,t)\|_{W^{1,q}(\Om)}
	+ \|\wab(\cdot,t)\|_{L^\infty(\Om)}
	\le C(q,K,\del)
	\qquad \mbox{for all } t\in (t_0,\tm).
  \ee
\end{lem}
\proof
  Given $n\le 3$ and $q>n$, we have $(n-3)q<n$ and hence $(n-2)q<n+q$, that is, $\frac{nq}{n+q}<\frac{n}{(n-2)_+}$.
  We can therefore fix $p=p(q)>1$ fulfilling
  \be{31.6}
	\frac{nq}{n+q} < p < \frac{n}{(n-2)_+},
  \ee
  and using that $q>n$, we can choose $r=r(q)>n$ in such a way that $r<q$.
  We may then rely on known smoothing properties of the Neumann heat semigroup $(e^{t\Del})_{t\ge 0}$ on $\Om$
  (\cite{win_JDE}, \cite{FIWY}) to find positive constants $\kappa$ and $c_i=c_i(q), i\in\{1,...,5\}$, such that whenever $t>0$,
  \be{31.77}
	\|e^{t\Del}\vp\|_{L^p(\Om)} \le c_1 \|\vp\|_{L^\infty(\Om)}
	\qquad \mbox{for all } \vp\in C^0(\bom)
  \ee
  and
  \be{31.7}
	\|e^{t\Del}\vp\|_{L^p(\Om)} \le c_1 \cdot \Big(1+t^{-\frac{n}{2}(1-\frac{1}{p})}\Big) \|\vp\|_{L^1(\Om)}
	\qquad \mbox{for all } \vp\in C^0(\bom)
  \ee
  as well as
  \be{31.8}
	\|e^{t\Del}\vp\|_{W^{1,q}(\Om)} \le c_3 \|\vp\|_{W^{1,q}(\Om)}
	\qquad \mbox{for all } \vp\in W^{1,q}(\Om)
  \ee
  and
  \be{31.9}
	\|e^{t\Del}\vp\|_{W^{1,q}(\Om)} \le c_4 \cdot \Big( 1+t^{-\frac{1}{2}-\frac{n}{2}(\frac{1}{p}-\frac{1}{q})} \Big) 
		\|\vp\|_{L^p(\Om)}
	\qquad \mbox{for all } \vp\in C^0(\bom)
  \ee
  and
  \be{31.10}
	\|e^{t\Del}\na\cdot \vp\|_{L^\infty(\Om)} \le c_5 \cdot \Big( 1+t^{-\frac{1}{2}-\frac{n}{2r}}\Big) e^{-\kappa t}
		\|\vp\|_{L^r(\Om)}
	\qquad \mbox{for all $\vp\in C^1(\bom;\R^n)$ such that $\vp\cdot\nu=0$ on $\pO$.}
  \ee
  Then assuming that $\del>0$ and $k>0$, that $\al>0$ and $\bet\ge\del$, and that (\ref{31.2}) and (\ref{31.3}) hold with some
  $t_0\in [0,\tm)$,
  we first employ a variation-of-constants representation associated with the third equation in (\ref{0ab}) to see that thanks
  to (\ref{31.77}), (\ref{31.7}), (\ref{31.2}), (\ref{31.3}) and the inequality $\bet\ge\del$,
  \bea{31.11}
	\|\wab(\cdot,t)\|_{L^p(\Om)}
	&=& \bigg\| e^{(t-t_0)(\Del-\bet)} \wab(\cdot,t_0)
	+ \int_{t_0}^t e^{(t-s)(\Del-\bet)} \big\{ \al \uab(\cdot,s)\big\} ds \bigg\|_{L^p(\Om)} \nn\\
	&\le& e^{-\bet (t-t_0)} \big\| e^{(t-t_0)\Del} \wab(\cdot,t_0)\big\|_{L^p(\Om)}
	+ \int_{t_0}^t e^{-\bet (t-s)} \big\| e^{(t-s)\Del} \big\{ \al\uab(\cdot,s)\big\} \big\|_{L^p(\Om)} ds \nn\\
	&\le& c_1 \|w(\cdot,t_0)\|_{L^\infty(\Om)}
	+ c_2 \int_{t_0}^t \Big(1+(t-s)^{-\frac{n}{2}(1-\frac{1}{p})} \Big) e^{-\bet (t-s)} \|\al\uab(\cdot,s)\|_{L^1(\Om)} ds \nn\\
	&\le& c_1 K 
	+ c_2 K \int_{t_0}^t \Big(1+(t-s)^{-\frac{n}{2}(1-\frac{1}{p})} \Big) e^{-\bet (t-s)} ds \nn\\
	&\le& c_6 \equiv c_6(q,K,\del):=
	c_1 K + c_2 K \int_0^\infty \Big(1+\sig^{-\frac{n}{2}(1-\frac{1}{p})}\big) e^{-\del\sig} d\sig
  \eea
  for all $t\in (t_0,\tm)$,
  where $c_6$ is finite, because $p>\frac{nq}{n+q}=\frac{n}{\frac{n}{q}+1}>\frac{n}{2}$ according to (\ref{31.6}) and the hypothesis
  $q>n$.\abs
  We next use this together with (\ref{31.8}), (\ref{31.9}) and again (\ref{31.3}) to find that
  for all $t\in (t_0,\tm)$,
  \bea{31.12}
	\|\vab(\cdot,t)\|_{W^{1,q}(\Om)}
	&=& \bigg\| e^{(t-t_0)(\Del-1)} \vab(\cdot,t_0) + \int_{t_0}^t e^{(t-s)(\Del-1)} \wab(\cdot,s) ds \bigg\|_{W^{1,q}(\Om)} 
		\nn\\
	&\le& e^{-(t-t_0)} \|e^{(t-t_0)\Del} \vab(\cdot,t_0)\|_{W^{1,q}(\Om)}
	+ \int_{t_0}^t e^{-(t-s)} \|e^{(t-s)\Del} \wab(\cdot,s)\|_{W^{1,q}(\Om)} ds \nn\\
	&\le& c_3 \|\vab(\cdot,t_0)\|_{W^{1,q}(\Om)}
	+ c_4 \int_{t_0}^t \Big(1+(t-s)^{-\frac{1}{2}-\frac{n}{2}(\frac{1}{p}-\frac{1}{q})}\Big) e^{-(t-s)} 
		\|\wab(\cdot,s)\|_{L^p(\Om)} ds \nn\\
	&\le& c_3 K 
	+ c_4 c_6 \int_{t_0}^t \Big(1+(t-s)^{-\frac{1}{2}-\frac{n}{2}(\frac{1}{p}-\frac{1}{q})}\Big) e^{-(t-s)} ds \nn\\
	&\le& c_7\equiv c_(q,K,\del)
		:=c_3 K + c_4 c_6 \int_0^\infty \Big(1+\sig^{-\frac{1}{2}-\frac{n}{2}(\frac{1}{p}-\frac{1}{q})} \Big) e^{-\sig} 
		d\sig,
  \eea
  with finiteness of $c_7$ being ensured by the left inequality in (\ref{31.6}), which precisely asserts that 
  $\frac{1}{2}+\frac{n}{2}(\frac{1}{p}-\frac{1}{q})
  <\frac{1}{2}+\frac{n}{2}(\frac{n+q}{nq}-\frac{1}{q}) =1$, namely.\abs
  Now (\ref{31.12}) in turn enables us to estimate the finite numbers
  \bas
	M(T):=\sup_{t\in (t_0,T)} \|\zab(\cdot,t)\|_{L^\infty(\Om)},
	\qquad T\in (t_0,\tm),
  \eas
  by drawing on a Duhamel representation of $\zab\equiv \uab+\wab$ corresponding to (\ref{0z}).
  Indeed, by means of the comparison principle and (\ref{31.10}), we obtain that
  for all $t\in (t_0,\tm)$,
  \bea{31.13}
	\|\zab(\cdot,t)\|_{L^\infty(\Om)}
	&=& \bigg\| e^{(t-t_0)\Del} \zab(\cdot,t_0) - \int_{t_0}^t e^{(t-s)\Del} \na\cdot \big\{ \uab(\cdot,s)\na \vab(\cdot,s)\big\}
		ds \bigg\|_{L^\infty(\Om)} \nn\\
	&\le& \|\zab(\cdot,t_0)\|_{L^\infty(\Om)} \nn\\
	& & + \int_{t_0}^t \Big(1+(t-s)^{-\frac{1}{2}-\frac{n}{2r}}\Big) e^{-\kappa(t-s)} 
		\big\| \uab(\cdot,s)\na\vab(\cdot,s)\|_{L^r(\Om)} ds,
  \eea
  where combining the H\"older inequality with (\ref{31.12}) and the ordering $0\le \uab\le \zab$ as well as the fact that
  $\|\zab(\cdot,t)\|_{L^1(\Om)}=\|\zab(\cdot,t_0)\|_{L^1(\Om)} \le |\Om| \cdot \|\zab(\cdot,t_0)\|_{L^\infty(\Om)}$
  for all $t\in (t_0,\tm)$ by (\ref{mass}),
  \bas
	& & \hs{-20mm}
	\big\|\uab(\cdot,s)\na\vab(\cdot,s)\big\|_{L^r(\Om)} \\
	&\le& \|\uab(\cdot,s)\|_{L^\frac{qr}{q-r}(\Om)} \|\na\vab(\cdot,s)\|_{L^q(\Om)} \\
	&\le& c_7 \|\uab(\cdot,s)\|_{L^\frac{qr}{q-r}(\Om)} \\
	&\le& c_7\|\zab(\cdot,s)\|_{L^\frac{qr}{q-r}(\Om)} \\
	&\le& c_7 \|\zab(\cdot,s)\|_{L^\infty(\Om)}^\theta \|\zab(\cdot,s)\|_{L^1(\Om)}^{1-\theta} \\
	&\le& c_7 |\Om|^{1-\theta} \|\zab(\cdot,t_0)\|_{L^\infty(\Om)}^{1-\theta} M^\theta(T)
	\qquad \mbox{for all $s\in (t_0,T)$ and } T\in (t_0,\tm),
  \eas
  where $\theta\equiv \theta(q):=1-\frac{q-r}{qr} \in (0,1)$.
  In light of (\ref{31.4}), from (\ref{31.13}) we hence infer that if we abbreviate 
  $c_8\equiv c_8(q,K,\del):=\max\big\{ K \, , \, 
  c_7 |\Om|^{1-\theta} K^{1-\theta} \int_0^\infty (1+\sig^{-\frac{1}{2}-\frac{n}{2r}}) e^{-\kappa\sig} d\sig \big\}$
  and note that $c_8$ is finite since $r>n$, then
  \bas
	\|\zab(\cdot,t)\|_{L^\infty(\Om)}
	&\le& K + c_7 |\Om|^{1-\theta} K^{1-\theta} M^\theta(T) 
		\int_{t_0}^t \Big(1+(t-s)^{-\frac{1}{2}-\frac{n}{2r}}\Big) e^{-\kappa(t-s)} ds \\
	&\le& c_8 +c_8 M^\theta(T)
	\qquad \mbox{for all $t\in (t_0,T)$ and any } T\in (t_0,\tm).
  \eas
  Therefore, 
  \bas
	M(T) \le c_8 + c_8 M^\theta(T)
	\qquad \mbox{for all } T\in (t_0,\tm),
  \eas
  which implies that
  \bas
	M(T) \le c_9\equiv c_9(q,K,\del):=\max \Big\{ 1 \, , \, (2c_8)^\frac{1}{1-\theta} \Big\}	
	\qquad \mbox{for all } T\in (t_0,\tm),
  \eas
  because $\theta<1$. As thus
  \bas
	\|\zab(\cdot,t)\|_{L^\infty(\Om)} \le c_9
	\qquad \mbox{for all } t\in (t_0,\tm),
  \eas
  in view of (\ref{31.12}) we conclude that indeed (\ref{31.4}) holds if we let $C(q,K,\del):=c_7+c_9$.
\qed
Our announced result on unconditional global solvability and on boundedness of each individual trajectory in (\ref{0ab})
thereby indeed reduces to a corollary:\abs
\proofc of Proposition \ref{prop1}.\quad
  We only need to take $\tm$ and $(\uab,\vab,\wab)$ as provided by Lemma \ref{lem_loc}, fix any $q>n$, and note that, thanks
  to (\ref{mass}) and (\ref{init}), the hypotheses (\ref{31.2}) and (\ref{31.3}) are satisfied with $t_0:=0$ and the finite number
  $K\equiv K(q):=\max \Big\{ \al \cdot \big\{ \io u_0 + \io w_0\big\} \, , \, 
  \|u_0\|_{L^\infty(\Om)} + \|v_0\|_{W^{1,q}(\Om)} + \|w_0\|_{L^\infty(\Om)} \Big\}$.
  Therefore, Lemma \ref{lem31} becomes applicable with $\del:=\bet$ so as to guarantee that in view of (\ref{ext}) and
  (\ref{31.4}) we necessarily must have $\tm=\infty$, and that thereupon (\ref{31.4}) implies (\ref{1.2}).
\qed
\subsection{$\bet$-independent estimates for bounded $\al$}
The strength of Lemma \ref{lem31} is underlined by its ability to directly imply not only the above, but also our main result on boundedness in (\ref{0ab}) for arbitrarily large $\bet$, under the assumption that $\al$ remains
within some finite interval:\abs
\proofc of Proposition \ref{prop41}.\quad
  As (\ref{mass}) guarantees that
  \bas
	\al\|\uab(\cdot,t)\|_{L^1(\Om)}
	\le \al^\star \cdot \bigg\{ \io u_0 + \io w_0 \bigg\}
	\qquad \mbox{for all $t>0$, $\al\in (0,\al^\star]$ and } \bet>0,
  \eas
  an application of Lemma \ref{lem31} to $t_0:=0$ immediately yields the claim.
\qed
\subsection{$\al$-independent estimates for bounded $\bet$}
In comparison to Proposition \ref{prop41}, our path toward Proposition \ref{prop42} needs to be slightly more sophisticated, essentially a consequence of the crucial boundedness feature (\ref{31.2}) not seemingly being evident for large values of $\al$. To appropriately handle the associated challenges, let us first record a simple but crucial observation concerned with a parameter-independent feature of the corresponding mass functional, in particular one that provides some favourable control throughout any time interval of the form $(\tau,\infty)$ with $\tau>0$:
\begin{lem}\label{lem44}
  Let $n\le 3$, and assume (\ref{init}).
  Then 
  \be{44.1}
	\al \io \uab(\cdot,t) \le \frac{m}{et} + m\bet
	\qquad \mbox{for all $t>0, \al>0$ and } \bet>0,
  \ee
  where $m:=\io u_0+\io w_0$.
\end{lem}
\proof
  Given $\al>0$ and $\bet>0$, we abbreviate $y(t):=\io \uab(\cdot,t), t\ge 0$, and then obtain from the first equation in
  (\ref{0ab}) and (\ref{mass}) that
  \bas
	y'(t)
	&=& - \al \io \uab + \bet \io \wab \\
	&=& - \al y(t) + \bet \io \wab \\
	&\le& - \al y(t) + m\bet
	\qquad \mbox{for all } t>0.
  \eas
  An ODE comparison argument thus shows that
  \bas
	y(t)
	&\le& y(0) e^{-\al t}
	+ \int_0^t e^{-\al (t-s)} \cdot m\bet ds \\
	&=& y(0) e^{-\al t}
	+ \frac{m\bet}{\al} \cdot (1-e^{-\al t}) \\
	&\le& me^{-\al t}
	+ \frac{m\bet}{\al}
	\qquad \mbox{for all } t>0,
  \eas
  because $y(0)\le m$. 
  Using that $\xi e^{-\xi} \le \frac{1}{e}$ for all $\xi\ge 0$, we hence infer that
  \bas
	\al y(t) 
	&\le& \frac{m\cdot \al t e^{-\al t}}{t} + m\bet \\
	&\le& \frac{m}{et} + m\bet
	\qquad \mbox{for all } t>0,
  \eas
  and conclude as intended.
\qed
In view of Lemma \ref{lem44}, the statement of Proposition \ref{prop42} will result from Lemma \ref{lem31} as soon as disadvantageous initial layer formation can be ruled out. This can indeed be achieved in the course of a self-map type regularity reasoning concerned with the subsystem of (\ref{0ab})-(\ref{0z}) solved by $(\zab,\vab)$, acting within suitably small but parameter-independent time intervals:
\begin{lem}\label{lem43}
  Let $n\le 3$, and assume (\ref{init}).
  Then for all $q>n$ there exist $T(q)>0$ and $C(q)>0$ such that
  \bea{43.1}
	& & \hs{-20mm}
	\|\uab(\cdot,t)\|_{L^\infty(\Om)}
	+ \|\vab(\cdot,t)\|_{W^{1,q}(\Om)}
	+ \|\wab(\cdot,t)\|_{L^\infty(\Om)}
	\le C(q,\al^\star,\del) \nn\\[1mm]
	& & \hs{40mm}
	\qquad \mbox{for all $t\in [0,T(q)]$ and any $\al>0$ and $\bet>0$.}
  \eea
\end{lem}
\proof
  We again employ standard smoothing estimates for the Neumann heat semigroup $(e^{t\Del})_{t\ge 0}$ to fix
  $c_i=c_i(q), i\in\{1,2,3\}$, with the property that if $t\in (0,1)$, then
  \be{43.2}	
	\|e^{t\Del}\vp\|_{W^{1,q}(\Om)} \le c_1 \|\vp\|_{W^{1,q}(\Om)}
	\qquad \mbox{for all } \vp\in W^{1,q}(\Om)
  \ee
  and
  \be{43.3}
	\|e^{t\Del}\vp\|_{W^{1,q}(\Om)} \le c_2 t^{-\frac{1}{2}} \|\vp\|_{L^\infty(\Om)}
	\qquad \mbox{for all } \vp\in C^0(\bom)
  \ee
  as well as
  \be{43.4}
	\|e^{t\Del}\na\cdot\vp\|_{L^\infty(\Om)}
	\le c_3 t^{-\frac{1}{2}-\frac{n}{2q}} \|\vp\|_{L^q(\Om)}
	\qquad \mbox{for all $\vp\in C^1(\bom;\R^n)$ such that $\vp\cdot\nu|_{\pO}=0$,}
  \ee
  and noting that $\frac{1}{2}-\frac{n}{2q}$ is positive, we thereupon choose $T=T(q)\in (0,1]$ small enough fulfilling
  \be{43.5}
	c_3 c_4 (c_5+1) \cdot \frac{T^{\frac{1}{2}-\frac{n}{2q}}}{\frac{1}{2}-\frac{n}{2q}} \le \frac{1}{2},
  \ee
  where
  \be{43.6}
	c_4\equiv c_4(q):=c_1 \|v_0\|_{W^{1,q}(\Om)} + 2c_2 \cdot (c_5+1)
  \ee
  with
  \be{43.7}
	c_5:=\|u_0+w_0\|_{L^\infty(\Om)}.
  \ee
  By continuity of the functions from (\ref{z}), for each $\al>0$ and $\bet>0$ the set
  \bas
	\sab:=\Big\{ T'\in (0,T) \ \Big| \ \|\zab(\cdot,t)\|_{L^\infty(\Om)} \le c_5+1 \mbox{ for all } t\in (0,T') \Big\}
  \eas
  then is not empty and hence $\tab:=\sup \sab$ is a well-defined element of $(0,T]$. To make sure that actually $\tab=T$,
  we first combine (\ref{43.2}) and (\ref{43.3}) with the inequalities $0\le\wab\le\zab$ and (\ref{43.6}) to see that
  \bea{43.8}
	\|\vab(\cdot,t)\|_{W^{1,q}(\Om)}
	&=& \bigg\| e^{t(\Del-1)} v_0 
	+ \int_0^t e^{(t-s)(\Del-1)} \wab(\cdot,s) ds \bigg\|_{W^{1,q}(\Om)} \nn\\
	&\le& e^{-t} \|e^{t\Del} v_0\|_{W^{1,q}(\Om)}
	+ \int_0^t e^{-(t-s)} \|e^{(t-s)\Del} \wab(\cdot,s)\|_{W^{1,q}(\Om)} ds \nn\\
	&\le& \|e^{t\Del} v_0\|_{W^{1,q}(\Om)}
	+ \int_0^t \|e^{(t-s)\Del} \wab(\cdot,s)\|_{W^{1,q}(\Om)} ds \nn\\
	&\le& c_1 \|v_0\|_{W^{1,q}(\Om)}
	+ c_2 \int_0^t (t-s)^{-\frac{1}{2}} \|\wab(\cdot,s)\|_{L^\infty(\Om)} ds \nn\\
	&\le& c_1 \|v_0\|_{W^{1,q}(\Om)}
	+ c_2 \int_0^t (t-s)^{-\frac{1}{2}} \|\zab(\cdot,s)\|_{L^\infty(\Om)} ds \nn\\
	&\le& c_1 \|v_0\|_{W^{1,q}(\Om)}
	+ c_2 (c_5+1) \int_0^t (t-s)^{-\frac{1}{2}} ds \nn\\
	&=& c_1 \|v_0\|_{W^{1,q}(\Om)}
	+ c_2 (c_5+1) \cdot \frac{t^\frac{1}{2}}{\frac{1}{2}} \nn\\
	&\le& c_4
	\qquad \mbox{for all } t\in (0,\tab),
  \eea
  because $\tab\le T \le 1$.
  Therefore, in view of the comparison principle, (\ref{43.4}), (\ref{43.6}) and (\ref{43.8}) and the fact that
  $0\le\uab\le\zab$, using (\ref{0z}) and (\ref{43.5}) we infer that
  \bas
	\|\zab(\cdot,t)\|_{L^\infty(\Om)}
	&=& \bigg\| e^{t\Del} (u_0+w_0)
	- \int_0^t e^{(t-s)\Del} \na\cdot\big\{ \uab(\cdot,s)\na\vab(\cdot,s)\big\} ds \bigg\|_{L^\infty(\Om)} \\
	&\le& \|e^{t\Del} (u_0+w_0)\|_{L^\infty(\Om)}
	+ \int_0^t \big\| e^{(t-s)\Del} \na\cdot\big\{ \uab(\cdot,s)\na\vab(\cdot,s)\big\} \big\|_{L^\infty(\Om)} ds \\
	&\le& \|u_0+w_0\|_{L^\infty(\Om)}
	+ c_3 \int_0^t (t-s)^{-\frac{1}{2}-\frac{n}{2q}} \big\|\uab(\cdot,s)\na\vab(\cdot,s)\big\|_{L^q(\Om)} ds \\
	&\le& \|u_0+w_0\|_{L^\infty(\Om)}
	+ c_3 \int_0^t (t-s)^{-\frac{1}{2}-\frac{n}{2q}} \|\uab(\cdot,s)\|_{L^\infty(\Om)} \|\na\vab(\cdot,s)\|_{L^q(\Om)} ds \\
	&\le& \|u_0+w_0\|_{L^\infty(\Om)}
	+ c_3 \int_0^t (t-s)^{-\frac{1}{2}-\frac{n}{2q}} \|\zab(\cdot,s)\|_{L^\infty(\Om)} \|\na\vab(\cdot,s)\|_{L^q(\Om)} ds \\
	&\le& \|u_0+w_0\|_{L^\infty(\Om)}
	+ c_3 c_4 (c_5+1) \int_0^t (t-s)^{-\frac{1}{2}-\frac{n}{2q}} ds \\
	&=& \|u_0+w_0\|_{L^\infty(\Om)}
	+ c_3 c_4 (c_5+1) \cdot \frac{t^{\frac{1}{2}-\frac{n}{2q}}}{\frac{1}{2}-\frac{n}{2q}} \\
	&\le& \|u_0+w_0\|_{L^\infty(\Om)}
	+ \frac{1}{2}
	\qquad \mbox{for all } t\in (0,\tab),
  \eas
  whence again by continuity of $\zab$, it follows that indeed $\tab$ cannot be smaller than $T$. Together with (\ref{43.8}),
  the definition of $\sab$ thus implies (\ref{43.1}) if we let $C(q):=c_4+c_5+1$.
\qed
Our main result on $\al$-independent boundedness in (\ref{0ab}) for parameters $\bet$ from fixed compact 
subintervals of $(0,\infty)$ can now be obtained on once more employing Lemma \ref{lem31}, this time on the 
basis of the previous two lemmata.\abs
\proofc of Proposition \ref{prop42}. \quad
  We let $T(q)>0$ be as provided by Lemma \ref{lem43}, and take any $t_0=t_0(q)\in (0,T(q))$.
  Then (\ref{43.1}) say that with some $c_1=c_1(q)>0$ we have
  \bea{42.2}
	& & \hs{-20mm}
	\|\uab(\cdot,t)\|_{L^\infty(\Om)}
	+ \|\vab(\cdot,t)\|_{W^{1,q}(\Om)}
	+ \|\wab(\cdot,t)\|_{L^\infty(\Om)}
	\le c_1 \nn\\[1mm]
	& & \hs{40mm}
	\qquad \mbox{for all $t\in [0,t_0]$ and any $\al>0$ and } \bet>0,
  \eea
  whereupon Lemma \ref{lem44} ensures that
  \be{42.3}
	\al\|\uab(\cdot,t)\|_{L^1(\Om)}
	\le c_2\equiv c_2(q,\bet^\star)
	:=\frac{m}{et_0} + m\bet^\star
	\qquad \mbox{for all $t\ge t_0, \al>0$ and } \bet\in (0,\bet^\star].
  \ee
  Particularly relying on the inequality in (\ref{42.2}) when evaluated at $t=t_0$, 
  we may therefore once again invoke Lemma \ref{lem31} to find $c_3=c_3(q,\bet^\star,\del)>0$ such that
  \bas
	& & \hs{-20mm}
	\|\uab(\cdot,t)\|_{L^\infty(\Om)}
	+ \|\vab(\cdot,t)\|_{W^{1,q}(\Om)}
	+ \|\wab(\cdot,t)\|_{L^\infty(\Om)}
	\le c_3 \nn\\[1mm]
	& & \hs{40mm}
	\qquad \mbox{for all $t>t_0, \al>0$ and } \bet\in [\del,\bet^\star].
  \eas
  Combined with (\ref{42.2}), this establishes (\ref{42.1}).
\qed
\mysection{Unlimited growth. Proof of Theorem \ref{theo12}}\label{sect3}
To next address the unboundedness phenomenon formulated in Theorem \ref{theo12},
let us from now on concentrate on the one-parameter sub-family (\ref{0})
of (\ref{0ab}). In order to avoid unnecessarily abundant notation, 
we let $\ug,\vg,\wg$ and $\zg$ be as introduced
in Proposition \ref{prop1} and (\ref{z}), respectively, with $\al=\bet=\gam$ for arbitrary $\gam>0$.\abs
The argument taken to detect unlimited growth will pursue an indirect strategy, which at its core aims to ensure that under suitable hypotheses on parameter-independent boundedness features in (\ref{0}),
a certain relationship to a classical two-component Keller-Segel system 
can be established.
To foreshadow the contradiction that will eventually be achieved, let us import from the literature the following consequence of well-known results on the occurrence of finite-time blow-up in such problems.
\begin{lem}\label{lem11}
  Let $n\ge 2$, $R>0$ and $\Om=B_R(0)\subset\R^n$, and suppose that $\chi>0$ and $a>0$.
  Then there exist $T>0$ as well as radially symmetric nonnegative functions
  $z_0\in C^2(\bom)$ and $v_0\in C^2(\bom)$ such that the problem
  \be{11.1}
	\left\{ \begin{array}{ll}	
	z_t = \Del z - \chi\na\cdot (z\na v),
	\qquad & x\in\Om, \ t\in (0,T), \\[1mm]
	v_t = \Del v - v + az,
	\qquad & x\in\Om, \ t\in (0,T), \\[1mm]
	\frac{\pa z}{\pa\nu}=\frac{\pa v}{\pa\nu}=0,
	\qquad & x\in\pO, \ t>0, \\[1mm]
	z(x,0)=z_0(x), \quad v(x,0)=v_0(x), 
	\qquad & x\in\Om,
	\end{array} \right.
  \ee
  does not admit any classical solution $(z,v)$ with
  \bas
	\lbal
	0 \le z \in C^0(\bom\times [0,T]) \cap C^{2,1}(\bom\times (0,T))
	\qquad \mbox{and} \\[1mm]
	0 \le v \in\bigcap_{q>n} C^0([0,T];W^{1,q}(\Om)) \cap C^{2,1}(\bom\times (0,T)).
	\ear
  \eas
\end{lem}
\proof
  In view of a corresponding uniqueness property within the indicated class (\cite{horstmann_win}), 
  this readily follows from known arguments revealing the occurrence of finite-time blow-up in (\ref{11.1})
  when either $n=2$ (\cite{herrero_velazquez}) or $n\ge 3$ (\cite{win_JMPA}).
\qed
\subsection{Basic implications of presupposed boundedness properties}
Our considerations in this regard will be launched by the following basic observation on how presupposed
bounds for $\wg$ imply regularity of the taxis gradients.
\begin{lem}\label{lem51}
  Let $n\le 3$, and suppose that with some $T>0$ and $S\subset (0,\infty)$ we have
  \be{51.1}
	\sup_{\gamma\in S} \sup_{t\in (0,T)} \|\wg(\cdot,t)\|_{L^\infty(\Om)} <\infty.
  \ee
  Then for all $q>n$ there exists $C(q)>0$ such that
  \be{51.2}
	\|\vg(\cdot,t)\|_{W^{1,q}(\Om)} \le C(q)
	\qquad \mbox{for all $t\in (0,T)$ and } \gam\in S.
  \ee
\end{lem}
\proof
  This can be seen by again relying on known regularization features of the Neumann heat semigroup $(e^{t\Del})_{t\ge 0}$
  on $\Om$, according to which, namely, we can find $c_1=c_1(q)>0$ and $c_2=c_2(q)>0$ such that
  \bas
	\|\vg(\cdot,t)\|_{W^{1,q}(\Om)}
	&\le& c_1 e^{-t} \|v_0\|_{W^{1,q}(\Om)}
	+ c_1\int_0^t \Big(1+(t-s)^{-\frac{1}{2}}\Big) e^{-(t-s)} \|\wg(\cdot,s)\|_{L^\infty(\Om)} ds \\
	&\le& c_2 + c_2\cdot\sup_{s\in (0,t)} \|\wg(\cdot,s)\|_{L^\infty(\Om)}
  \eas
  for all $t>0$ and $\gam>0$.
\qed
In addition, a standard testing procedure shows that
the hypothesis in (\ref{51.1}) also entails space-time $L^2$ regularity features of $\Del\vg$ and $v_{\gam t}$.
\begin{lem}\label{lem511}
  Let $n\le 3$, and assume (\ref{51.1}) with some $T>0$ and $S\subset (0,\infty)$.
  Then one can fix $C>0$ in such a way that
  \be{511.1}
	\int_0^T \io |\Del\vg|^2 \le C
	\qquad \mbox{for all } \gam\in S
  \ee
  and
  \be{511.2}
	\int_0^T \io v_{\gam t}^2 \le C
	\qquad \mbox{for all } \gam\in S.
  \ee
\end{lem}
\proof
  We test the second equation in (\ref{0}) by $-\Del\vg$ and use Young's inequality as well as (\ref{51.1}) to find that with some
  $c_1>0$,
  \bas
	\frac{1}{2} \frac{d}{dt} \io |\na\vg|^2
	+ \io |\Del\vg|^2
	&=& - \io |\na\vg|^2
	+ \io \wg\Del\vg \\
	&\le& \frac{1}{2} \io |\Del\vg|^2
	+ \frac{1}{2} \io \wg^2 \\
	&\le& \frac{1}{2} \io |\Del\vg|^2
	+ c_1
	\qquad \mbox{for all $t\in (0,T)$ and } \gam\in S,
  \eas
  and that hence
  \be{511.3}
	\io \big|\na\vg(\cdot,T)\big|^2 
	+ \int_0^T \io |\Del\vg|^2
	\le \io |\na v_0|^2
	+ 2c_1 T
	\qquad \mbox{for all } \gam\in S.
  \ee
  As from Lemma \ref{lem51} and the hypothesis (\ref{51.1}) we furthermore obtain $c_2>0$ satisfying
  \bas
	v_{\gam t}^2
	= (\Del\vg-\vg+\wg)^2
	\le 3|\Del\vg|^2
	+3\vg^2 + 3\wg^2
	\le 3|\Del\vg|^2 +c_2
	\quad \mbox{in } \Om\times (0,T)
	\qquad \mbox{for all } \gam\in S,
  \eas
  the inequality in (\ref{511.3}) implies both (\ref{511.1}) and (\ref{511.2}).
\qed
Now once more thanks to smoothing features of the Neumann heat semigroup, Lemma \ref{lem51} can be seen to imply
uniform $L^\infty$ bounds for $\ug+\wg$ under the assumption in (\ref{51.1}):
\begin{lem}\label{lem52}
  Given $n\le 3$ and assuming (\ref{51.1}) to hold with some $T>0$ and $S\subset (0,\infty)$,
  one can find $C>0$ such that if for $\gam>0$ we let $\zg\equiv \ug+\wg$ be as defined in (\ref{z}) with $\al=\bet=\gam$, then
  \be{52.1}
	\|\zg(\cdot,t)\|_{L^\infty(\Om)} \le C
	\qquad \mbox{for all $t\in (0,T)$ and } \gam\in S.
  \ee
\end{lem}
\proof
  We pick any $q>n$ and then take $r\in (n,q)$, and similar to the reasoning in Lemma \ref{lem31}, we introduce
  \bas
	M_\gam(T'):=\sup_{t\in (0,T')} \|\zg(\cdot,t)\|_{L^\infty(\Om)},
	\qquad T'\in (0,T), \ \gam\in S.
  \eas
  To estimate these quantities, we once more draw on smoothing properties of the Neumann heat semigroup $(e^{t\Del})_{t\ge 0}$ 
  on $\Om$, and on Lemma \ref{lem51} and (\ref{mass}) to fix positive constants $\kappa, c_1$ and $c_2$ such that 
  whenever $\gam\in S$,
  \bas
	\|\zg(\cdot,t)\|_{L^\infty(\Om)}
	&\le& \big\| e^{t\Del}(u_0+w_0)\big\|_{L^\infty(\Om)}
	+ \int_0^t \big\| e^{(t-s)\Del} \na\cdot \big\{ \ug(\cdot,s)\na \vg(\cdot,s)\big\} \big\|_{L^\infty(\Om)} ds \\
	&\le& \|u_0+w_0\|_{L^\infty(\Om)}
	+ c_1\int_0^t \Big(1+(t-s)^{-\frac{1}{2}-\frac{n}{2r}}\Big) e^{-\kappa(t-s)} \|\ug(\cdot,s)\na\vg(\cdot,s)\|_{L^r(\Om)} ds \\
	&\le& \|u_0+w_0\|_{L^\infty(\Om)} \\
	& & + c_1\int_0^t \Big(1+(t-s)^{-\frac{1}{2}-\frac{n}{2r}}\Big) e^{-\kappa(t-s)} \|\zg(\cdot,s)\|_{L^\infty(\Om)}^\theta 
		\|\zg(\cdot,s)\|_{L^1(\Om)}^{1-\theta} \|\na\vg(\cdot,s)\|_{L^q(\Om)} ds \\
	&\le& \|u_0+w_0\|_{L^\infty(\Om)} \\
	& & + c_1 \|u_0+w_0\|_{L^1(\Om)}^{1-\theta} \cdot \bigg\{ \sup_{s\in (0,T)} \|\na\vg(\cdot,s)\|_{L^q(\Om)} \bigg\} \cdot
		M_\gam^\theta(T') \times \\
	& & \hs{60mm} 
	\times \int_0^t \Big(1+(t-s)^{-\frac{1}{2}-\frac{n}{2r}}\Big) e^{-\kappa(t-s)} ds \\
	&\le& c_2 + c_2 M_\gam^\theta(T')
	\qquad \mbox{for all $t\in (0,T')$ and } T'\in (0,T),
  \eas
  where $\theta:=1-\frac{q-r}{qr}$. Since $\theta\in (0,1)$, the inequality
  $M_\gam(T') \le c_2 + c_2 M_\gam^\theta(T')$, as thereby implied for each $T'\in (0,T)$ and $\gam\in S$, entails that
  $M_\gam(T') \le \max\big\{ 1 \, , \, (2c_2)^\frac{1}{1-\theta}\big\}$ for any such $T'$ and $\gam$, and hence establishes
  the claim.
\qed
On the basis of the parabolic problems solved by $\zg$ and $\vg$, the estimates from Lemma \ref{lem51} and Lemma \ref{lem52}
quite directly imply bound also in H\"older spaces.
\begin{lem}\label{lem53}
  If $n\le 3$ and (\ref{51.1}) is valid with some $T>0$ and $S\subset (0,\infty)$,
  then there exist $\vartheta\in (0,1)$ and $C>0$ such that
  \be{53.1}
	\|\vg(\cdot,t)\|_{C^{\vartheta,\frac{\vartheta}{2}}(\bom\times [0,T])} \le C
	\qquad \mbox{for all } \gam\in S
  \ee
  and
  \be{53.2}
	\|\zg(\cdot,t)\|_{C^{\vartheta,\frac{\vartheta}{2}}(\bom\times [0,T])} \le C
	\qquad \mbox{for all } \gam\in S,
  \ee
  where $(\zg)_{\gam>0}$ is as in Lemma \ref{lem52}.
\end{lem}
\proof
  This readily follows from standard results on H\"older regularity in scalar parabolic equations (\cite{porzio_vespri}),
  based on the hypothesis in (\ref{51.1}) in the derivation of (\ref{53.1}) and on the outcomes of Lemma \ref{lem51} and 
  Lemma \ref{lem52}, the former being applied to any $q>n$, in the verification of (\ref{53.2}).
\qed
\subsection{Controlling the difference $\wg-\ug$ for large $\gam$}
Forming the cornerstone of our reasoning related to Theorem \ref{theo12}, this section will reveal that
under the assumption that (\ref{51.1}) be valid for some unbounded set $S$, the difference $\wg-\ug$ necessarily
needs to be conveniently small for large values of $\gam$.
In contrast to those taken in most of what precedes, our arguments here will be more or less exclusively of variational nature.\abs
The first statement in this context will draw on a rather straightforward $L^2$ testing procedure.
\begin{lem}\label{lem54}
  Assume that $n\le 3$ and that $T>0$ and $S\subset (0,\infty)$ are such that (\ref{51.1}) holds.
  Then there exists $C>0$ such that the family $(\zg)_{\gam>0}$ from Lemma \ref{lem52} satisfies
  \be{54.1}
	\int_0^T \io |\na\zg|^2 \le C
	\qquad \mbox{for all } \gam\in S.
  \ee
\end{lem}
\proof
  On testing the respective version of (\ref{0z}) against $\zg$ in a straightforward manner, by means of Young's inequality we see 
  that
  \bas
	\frac{d}{dt} \io \zg^2 + 2 \io |\na\zg|^2
	&=& 2 \io \ug \na\vg\cdot\na\zg \\
	&\le& \io |\na\zg|^2
	+ \io \ug^2 |\na\vg|^2 \\
	&\le& \io |\na\zg|^2
	+ \|\zg\|_{L^\infty(\Om)}^2 \io |\na\vg|^2
	\qquad \mbox{for all $t>0$ and } \gam>0.
  \eas
  In view of Lemma \ref{lem51} and Lemma \ref{lem52}, we thus obtain that with some $c_1>0$,
  \bas
	\frac{d}{dt} \io \zg^2 + \io |\na \zg|^2 \le c_1
	\qquad \mbox{for all $t\in (0,T)$ and } \gam\in S,
  \eas
  and that thus
  \bas
	\io \zg^2(\cdot,T) + \int_0^T \io |\na\zg|^2
	\le \io (u_0+w_0)^2 + c_1 T
	\qquad \mbox{for all } \gam\in S,
  \eas
  from which (\ref{54.1}) follows.
\qed
A key step will now be prepared by a further quite elementary testing argument.
\begin{lem}\label{lem55}
  Let $n\le 3$, and given $\gam>0$, let $\zg$ be as in Lemma \ref{lem52}. Then
  \bea{55.1}
	& & \hs{-20mm}
	\frac{d}{dt} \bigg\{ \io \wg^2 - \io \wg \zg \bigg\}
	+ 2\io |\na\wg|^2 + \gam \io (\wg-\ug)^2 \nn\\
	&=& 2\io \na\wg\cdot\na\zg
	- \io \ug\na\vg\cdot\na\wg
	\qquad \mbox{for all } t>0.
  \eea
\end{lem}
\proof
  We multiply the third equation in (\ref{0}) by $\wg-\ug$ to see that
  \be{55.2}
	\io w_{\gam t} \cdot (\wg-\ug)
	= \io \Del\wg \cdot (\wg-\ug)
	- \gam \io (\wg-\ug)^2
	\qquad \mbox{for all } t>0,
  \ee
  where since $\wg-\ug\equiv 2\wg-\zg$,
  \be{55.3}
	\io \Del\wg \cdot (\wg-\ug)
	= - 2\io |\na\wg|^2
	+ \io \na\wg\cdot\na\zg
	\qquad \mbox{for all } t>0
  \ee
  and
  \bea{55.4}
	\io w_{\gam t} \cdot (\wg-\ug)
	&=& 2 \io w_{\gam t} \wg
	- \io w_{\gam t} \zg \nn\\
	&=& \frac{d}{dt} \io \wg^2
	- \frac{d}{dt} \io \wg \zg
	+ \io \wg z_{\gam t}
	\qquad \mbox{for all } t>0.
  \eea
  Since (\ref{0z}) implies that here 
  \bas
	\io \wg z_{\gam t}
	&=& \io \wg \Del\zg
	- \io \wg \na\cdot (\ug\na\vg) \\
	&=& - \io \na\wg\cdot\na\zg
	+ \io \ug\na\vg\cdot\na\wg
	\qquad \mbox{for all } t>0,
  \eas
  a combination of (\ref{55.2})-(\ref{55.4}) yields (\ref{55.1}).
\qed
The following consequence of the latter relies on the appearance of the potentially large factor $\gam$ on the left-hand side
of (\ref{55.1}).
\begin{lem}\label{lem56}
  Let $n\le 3$, and assume that $T>0$ and the unbounded set $S\subset (0,\infty)$ are such that (\ref{51.1}) holds.
  Then
  \be{56.1}
	\int_0^T \io (\wg-\ug)^2 \to 0
	\qquad \mbox{as } S\ni \gam \to\infty.
  \ee
\end{lem}
\proof
  According to the outcomes of Lemma \ref{lem51} and Lemma \ref{lem52} and our definition of $(\zg)_{\gam>0}$, we can fix
  $c_1>0, c_2>0$ and $c_3>0$ such that
  \be{56.2}
	\|\ug\|_{L^\infty(\Om)} \le c_1,
	\quad	
	\|\na\vg\|_{L^2(\Om)} \le c_2
	\quad \mbox{and} \quad
	\|\zg\|_{L^\infty(\Om)} \le c_3
	\qquad \mbox{for all $t\in [0,T]$ and } \gam\in S,
  \ee
  while Lemma \ref{lem54} provides $c_4>0$ fulfilling
  \be{56.3}
	\int_0^T \io |\na\zg|^2 \le c_4
	\qquad \mbox{for all } \gam\in S.
  \ee
  To make use of this in the context of (\ref{55.1}), we first employ Young's inequality to see that
  \bas
	2\io \na\wg\cdot\na\zg \le \io |\na\wg|^2 + \io |\na\zg|^2
	\qquad \mbox{for all $t>0$ and } \gam>0,
  \eas
  and that
  \bas
	-\io \ug\na\vg\cdot\na\wg
	&\le& \io |\na\wg|^2
	+ \io \ug^2 |\na\vg|^2 \\
	&\le& \io |\na\wg|^2
	+ \|\ug\|_{L^\infty(\Om)}^2 \io |\na\vg|^2
	\qquad \mbox{for all $t>0$ and } \gam>0.
  \eas
  An integration of (\ref{55.1}) therefore shows that again by Young's inequality, and by (\ref{56.2}) and (\ref{56.3}),
  \bas
	\gam \int_0^T \io (\wg-\ug)^2
	&\le& \io w_0^2 - \io w_0(u_0+w_0) \\
	& & - \io \wg^2(\cdot,T) + \io \wg(\cdot,T) \zg(\cdot,T) \\
	& & + \int_0^T \io |\na\zg|^2
	+ \int_0^T \|\ug(\cdot,t)\|_{L^\infty(\Om)}^2 \cdot \bigg\{ \io |\na\vg(\cdot,t)|^2 \bigg\} \, dt \\
	&\le& \io w_0^2 - \io w_0(u_0+w_0) \\
	& & +\frac{1}{4} \io \zg^2(\cdot,T) \\
	& & + \int_0^T \io |\na\zg|^2
	+ \int_0^T \|\ug(\cdot,t)\|_{L^\infty(\Om)}^2 \cdot \bigg\{ \io |\na\vg(\cdot,t)|^2 \bigg\} \, dt \\
	&\le& \io w_0^2 - \io w_0(u_0+w_0) \\
	& & +\frac{c_3^2 |\Om|}{4} \\
	& & + c_4 + c_1^2 c_2^2 T
	\qquad \mbox{for all } \gam\in S,
  \eas
  which directly results in (\ref{56.1}).
\qed
\subsection{A link to a classical Keller-Segel system. Conclusion}
On the basis of the estimates of the previous two sections, a straightforward subsequence extraction enables
us to construct a limit pair that forms a weak solution of a two-component Keller-Segel system, provided that (\ref{51.1})
holds for some unbounded $S$.
\begin{lem}\label{lem57}
  If $n\le 3$, and if (\ref{51.1}) is satisfied with some $T>0$ and some unbounded set $S\subset (0,\infty)$,
  then there exist $(\gam_j)_{j\in\N}\subset S$ as well as nonnegative functions
  \be{57.1}
	\lbal
	z_\infty \in \bigcup_{\vartheta\in (0,1)} C^{\vartheta,\frac{\vartheta}{2}}(\bom\times [0,T])
		\cap L^2((0,T);W^{1,2}(\Om))
	\qquad \mbox{and} \\[1mm]
	v_\infty \in \bigcup_{\vartheta\in (0,1)} C^{\vartheta,\frac{\vartheta}{2}}(\bom\times [0,T])
		\cap \bigcap_{q>n} L^\infty((0,T);W^{1,q}(\Om))
	\ear
  \ee
  such that $\gam_j\to\infty$ as $j\to\infty$, that with $(z_\gam)_{\gam>0}$ from Lemma \ref{lem52} we have
  \begin{eqnarray}
	\ug \wsto \frac{1}{2} z_\infty
	\qquad \mbox{in } L^\infty(\Om\times (0,T)),
	\label{57.2} \\
	\vg \to v_\infty
	\qquad \mbox{in } C^0(\bom\times [0,T]),
	\label{57.3} \\
	\na\vg \to \na v_\infty
	\qquad \mbox{in } L^2(\Om\times (0,T)),
	\label{57.4} \\
	\wg \wsto \frac{1}{2} z_\infty
	\qquad \mbox{in } L^\infty(\Om\times (0,T)),
	\label{57.5} \\
	\zg \to z_\infty
	\qquad \mbox{in } C^0(\bom\times [0,T])
	\qquad \qquad \mbox{and}
	\label{57.6} \\
	\na\zg \to \na z_\infty
	\qquad \mbox{in } L^2(\Om\times (0,T)),
	\label{57.7}
  \end{eqnarray}
  as $\gam=\gam_j\to\infty$, and that for each $\vp\in C_0^\infty(\bom\times [0,T))$,
  \be{wz}
	- \int_0^T \io z_\infty \vp_t
	- \io (u_0+v_0)\vp(\cdot,0)
	= - \int_0^T \io \na z_\infty\cdot\na\vp
	+ \frac{1}{2} \int_0^T \io z_\infty \na v_\infty\cdot\na\vp
  \ee
  and
  \be{wv}
	- \int_0^T \io v_\infty \vp_t
	- \io v_0 \vp(\cdot,0)
	= - \int_0^T \io \na v_\infty \cdot\na\vp
	- \int_0^T \io v_\infty \vp
	+ \frac{1}{2} \int_0^T \io z_\infty \vp.
  \ee
\end{lem}
\proof
  From Lemma \ref{lem51}, Lemma \ref{lem511} and Lemma \ref{lem53} we know that
  \bas
	(\vg)_{\gam\in S} 
	\mbox{ is bounded in $L^2((0,T);W^{1,2}(\Om))$ and in $L^\infty((0,T);W^{1,q}(\Om))$ for all $q>n$,}
  \eas
  and that moreover
  \bas
	(\vg)_{\gam\in S}
	\mbox{ is bounded in $C^{\vartheta_1,\frac{\vartheta_1}{2}}(\bom\times [0,T])$ for some $\vartheta_1\in (0,1)$,}
  \eas
  and that
  \bas
	(v_{\gam t})_{\gam\in S}
	\mbox{ is bounded in $L^2(\Om\times (0,T))$.}
  \eas
  Apart from that, Lemma \ref{lem53} together with Lemma \ref{lem54} ensures that
  \bas
	(\zg)_{\gam\in S}
	\mbox{ is bounded in $C^{\vartheta_2,\frac{\vartheta_2}{2}}(\bom\times [0,T])$ for some $\vartheta_2\in (0,1)$,}
  \eas
  and that
  \bas
	(\zg)_{\gam\in S} 
	\mbox{ is bounded in $L^2((0,T);W^{1,2}(\Om))$,}
  \eas
  while Lemma \ref{lem52} implies that both
  \bas
	(\ug)_{\gam\in S} 
	\mbox{ is bounded in $L^\infty(\Om\times (0,T))$}
  \eas
  and
  \bas
	(\wg)_{\gam\in S} 
	\mbox{ is bounded in $L^\infty(\Om\times (0,T))$}.
  \eas
  A straightforward extraction procedure that relies on the Aubin-Lions lemma (\cite{temam}), the Arzel\`a-Ascoli theorem and the
  unboundedness of $S$ thus yields $(\gam_j)_{j\in\N} \subset S$ such that $\gam_j\to\infty$ as $j\to\infty$,
  and that with some nonnegative functions $u_\infty, v_\infty, w_\infty$ and $z_\infty$
  fulfilling (\ref{57.1}) as well as $u_\infty \in L^\infty(\Om\times (0,T))$ and $w_\infty\in L^\infty(\Om\times (0,T))$,
  besides (\ref{57.3}), (\ref{57.4}), (\ref{57.6}) and (\ref{57.7}) we have
  \be{57.8}
	\ug \wsto u_\infty
	\quad \mbox{and} \quad
	\wg \wsto w_\infty
	\qquad \mbox{in } L^\infty(\Om\times (0,T))
  \ee
  as $\gam=\gam_j\to\infty$. 
  Since thus, in particular, $\wg-\ug \wto w_\infty-u_\infty$ in $L^2(\Om\times (0,T))$
  as $\gam=\gam_j\to\infty$,
  and since on the other hand from Lemma \ref{lem56} we know that $\wg-\ug\to 0$ in $L^2(\Om\times (0,T))$
  as $S\ni\gam\to \infty$, we obtain that necessarily
  \be{57.9}
	w_\infty=u_\infty
	\qquad \mbox{a.e.~in } \Om\times (0,T).
  \ee
  Independently, a combination of (\ref{57.6}) with the definition of $(\zg)_{\gam>0}$ shows that
  $\ug+\wg=\zg \wsto z_\infty$ in $L^\infty(\Om\times (0,T))$ as $\gam=\gam_j\to\infty$, and that therefore (\ref{57.8})
  implies that $u_\infty+w_\infty=z_\infty$ a.e.~in $\Om\times (0,T)$.
  From (\ref{57.9}) it thus follows that $u_\infty=w_\infty=\frac{1}{2} z_\infty$ a.e.~in $\Om\times (0,T)$, and that hence,
  due to (\ref{57.8}), also (\ref{57.2}) and (\ref{57.5}) hold
  as $\gam=\gam_j\to\infty$. \abs
  A derivation of the identities in (\ref{wz}) and (\ref{wv}) can therefore be achieved on the basis of (\ref{57.2})-(\ref{57.7})
  and the corresponding weak formulations associated with (\ref{0z}) and the second equation in (\ref{0}) in a standard manner.
\qed
In fact, any such limit pair solves the considered chemotaxis system even in the classical sense:
\begin{lem}\label{lem58}
  Let $n\le 3$, and suppose that (\ref{51.1}) holds with some $T>0$ and some unbounded set $S\subset (0,\infty)$,
  Then the functions obtained in Lemma \ref{lem57} have the additional properties that
  \be{58.1}
	\lbal
	z_\infty \in C^0(\bom\times [0,T]) \cap C^{2,1}(\bom\times (0,T))
	\qquad \mbox{and} \\[1mm]
	v_\infty \in \bigcap_{q>n} C^0([0,T];W^{1,q}(\Om)) \cap C^{2,1}(\bom\times (0,T)),
	\ear
  \ee
  and $(v_\infty,z_\infty)$ solves
  \be{58.2}
	\left\{ \begin{array}{ll}	
	\pa_t z_\infty= \Del z_\infty - \frac{1}{2} \na\cdot (z_\infty \na v_\infty),
	\qquad & x\in\Om, \ t\in (0,T), \\[1mm]
	\pa_t v_\infty = \Del v_\infty - v_\infty + \frac{1}{2} z_\infty,
	\qquad & x\in\Om, \ t\in (0,T), \\[1mm]
	\frac{\pa z_\infty}{\pa\nu}=\frac{\pa v_\infty}{\pa\nu}=0,
	\qquad & x\in\pO, \ t>0, \\[1mm]
	z_\infty(x,0)=u_0(x)+w_0(x), \quad v_\infty(x,0)=v_0(x), 
	\qquad & x\in\Om,
	\end{array} \right.
  \ee
  in the classical sense.
\end{lem}
\proof
  As Lemma \ref{lem57} ensures that $v_\infty$ solves the second sub-problem in (\ref{58.2}) in the standard weak sense addressed,
  e.g., in \cite{LSU}, the H\"older continuity property of $z_\infty$ asserted by (\ref{57.1}) enables us to employ classical
  parabolic regularity theory of Schauder type (\cite{LSU}) to conclude that $v_\infty$ indeed has the smoothness features
  in (\ref{58.1}) and solves its respective part in (\ref{58.2}) in the claimed pointwise sense.
  The corresponding properties of $z_\infty$ can thereupon be verified in a quite similar manner.
\qed
As a final ingredient to our proof of Theorem \ref{theo12}, let us perform a simple comparison argument in order to
make sure that all the above statements continue to hold if, instead of (\ref{51.1}), a corresponding bound on the functions
$\ug$ is assumed:
\begin{lem}\label{lem61}
  Let $n\le 3$, and assume that $T>0$ and $S\subset (0,\infty)$ are such that
  \be{61.1}
	\sup_{\gamma\in S} \sup_{t\in (0,T)} \|\ug(\cdot,t)\|_{L^\infty(\Om)} <\infty.
  \ee
  Then (\ref{51.1}) holds.
\end{lem}
\proof
  In line with (\ref{61.1}), we let $c_1>0$ be such that $\ug\le c_1$ in $\Om\times (0,T)$ for all $\gam\in S$, and writing
  $c_2:=\|w_0\|_{L^\infty(\Om)}$ we define $\ow(x,t):=\max\{ c_1 \, , \, c_2\}$ for $x\in\bom$ and $t\ge 0$.
  Then $\ow(x,0)\ge c_2 \ge \|w_0\|_{L^\infty(\Om)} \ge \wg(x,0)$ for all $x\in\Om$, and furthermore the inequality
  $\ow\ge c_1$ ensures that
  \bas
	\ow_t - \Del \ow + \gam \ow - \gam \ug
	&\gam \ow - \gam \ug
	\ge \gam c_1 - \gam \ug \ge 0
	\qquad \mbox{in } \Om\times (0,T)
  \eas
  according to our choice of $c_1$. As clearly $\frac{\pa\ow}{\pa\nu}=0$ on $\pO\times (0,\infty)$, the comparison principle
  guarantees that $\wg \le \ow$ in $\Om\times (0,T)$ for all $\gam\in S$, and hence confirms (\ref{51.1}).
\qed
We can thereby readily derive our main result on the spontaneous emergence of arbitrarily large population densities
in (\ref{0}).\abs
\proofc of Theorem \ref{theo12}.
  We let $z_0\in W^{1,\infty}(\Om)$, $v_0\in W^{1,\infty}(\Om)$ and $T>0$ be as provided by Lemma \ref{lem11} when applied
  to $\chi:=\frac{1}{2}$ and $a:=\frac{1}{2}$, and let $u_0:=\frac{1}{2} z_0$ as well as $w_0:=\frac{1}{2} z_0$, for instance.
  Then given any unbounded $(\gam_j)_{j\in\N}\subset (0,\infty)$, invoking Lemma \ref{lem57} and Lemma \ref{lem58} with
  $S\equiv S^{(1)}:=\{\gam_j \ | \ j\in\N\}$ we obtain that due to Lemma \ref{lem11},
  the set $\big(\|w_{\gam_j}\|_{L^\infty(\Om\times (0,T))}\big)_{j\in\N}$ cannot be bounded, and that hence there must exist a
  subsequence $\big(\gam_{j_l^{(1)}}\big)_{l\in\N}$ of $(\gam_j)_{j\in\N}$ such that 
  \be{12.3}
	\big\|w_{\gam_{j_l^{(1)}}}\big\|_{L^\infty(\Om\times (0,T))} \to \infty
	\qquad \mbox{as } l\to\infty.
  \ee
  We thereupon apply Lemma \ref{lem61} and then again Lemmata \ref{lem57}, \ref{lem58} and \ref{lem11}, now
  to $S\equiv S^{(2)}:=\big\{\gam_{j_l^{(1)}} \ \big| \ l\in\N\big\}$, to infer that also
  $\Big( \big\| u_{\gam_{j_l^{(1)}}} \big\|_{L^\infty(\Om\times (0,T))} \Big)_{l\in\N}$ 
  must be unbounded. 
  We can thus extract a further subsequence $(\gam{j_k})_{k\in\N}$ of $\big(\gam_{j_l^{(1)}}\big)_{l\in\N}$ fulfilling (\ref{12.1}),
  and conclude by noting that then (\ref{12.3}) particularly entails (\ref{12.2}).
\qed

\bigskip

\bigskip

{\bf Acknowledgement.} \
  
  The first author acknowledges ``MIUR-Dipartimento di Eccellenza'' funding to the Dipartimento Interateneo di Scienze, Progetto e Politiche del Territorio (DIST).
  The second author acknowledges support of the {\em Deutsche Forschungsgemeinschaft} (Project No.~462888149).
\small

\end{document}